\documentclass[12pt]{amsart}
\usepackage{graphicx}
\usepackage{amssymb}
\usepackage{amsmath}
\usepackage{amsthm,amsfonts,bbm}
\usepackage{amscd}
\usepackage{geometry}
\usepackage[all,2cell]{xy}
\usepackage{epsfig,epstopdf}
\usepackage[colorlinks=true,anchorcolor=blue,filecolor=blue,linkcolor=blue,urlcolor=blue,citecolor=blue]{hyperref}
\usepackage{tikz}
\usepackage{caption}
\usepackage{subcaption}

\UseAllTwocells \SilentMatrices
\newtheorem{thm}{Theorem}[section]
\newtheorem{coro}[thm]{Corollary}
\newtheorem{lem}[thm]{Lemma}

\newtheorem{conj}[thm]{Conjecture}
\theoremstyle{definition}
\newtheorem{defi}[thm]{Definition}
\theoremstyle{remark}

\numberwithin{equation}{section}
\numberwithin{figure}{section}
\def\A{\mathcal{A}}
\def \am{\mathrm{am}}
\def\B{\mathcal{B}}
\def\C{\mathbb{C}}
\def\D{\mathcal{D}}
\def\diag{\mathrm{diag}}
\def \dim{\mathrm{dim}}
\def \F{\mathbb{F}}
\def \gm{\mathrm{gm}}
\def \i {\mathbf{i}}
\def \I{\mathcal{I}}
\def \K{\mathbb{K}}
\def \L{\mathcal{L}}
\def\la{\lambda}
\def \N{\mathbb{N}}
\def \p{\mathfrak{p}}
\def \part{\partial}
\def \PV{\mathbb{V}}
\def \Q{\mathcal{Q}}
\def \Res{\mathrm{Res}}
\def \S{\mathbb{S}}
\def \supp{\mathrm{supp}}
\def \spec{\mathrm{Spec}}
\def \Tr{\mathrm{Tr}}
\def \V{\mathcal{V}}
\def \Z{\mathbb{Z}}

\begin{document}
\title[Multiplicity of eigenvalues]{The multiplicity of eigenvalues of nonnegative tensors and hypergraphs}

\author[Y.-Z. Fan]{Yi-Zheng Fan*}
\address{Center for Pure Mathematics, School of Mathematical Sciences, Anhui University, Hefei 230601, P. R. China}
\email{fanyz@ahu.edu.cn}
\thanks{*Corresponding author.
This work was supported by the National Natural Science Foundation of China (Grant No. 12331012).}



\begin{abstract}
Hu and Ye conjectured that for an  $n$-dimensional tensor $\A$ of order $k$ with an eigenvalue $\la$ and the corresponding eigenvariety $\V_\la(\A)$, the algebraic multiplicity  $\am(\la)$ of $\la$ satisfies:
$$\am(\la) \ge \sum_{i=1}^\kappa \dim(V_i)(k-1)^{\dim(V_i)-1},$$
where $V_1,\ldots,V_\kappa$ are all irreducible components of $\V_\la(\A)$.
In this paper, we establish that for any weakly irreducible nonnegative tensor $\A$ with spectral radius $\rho$, all the eigenvalues $\la$ of $\A$ with $|\la|=\rho$ satisfy $\am(\la) \ge |\PV_\la(\A)|$, where $\PV_\la(\A)$ is the projective eigenvariety of $\A$ associated with $\la$.
As a direct consequence, we confirm the Hu-Ye Conjecture for two classes of eigenvalues: (1) all  eigenvalues of  weakly irreducible nonnegative tensors with modulus equal to the spectral radius, and (2)  the least H-eigenvalues of weakly irreducible $Z$-tensors.
Furthermore, we characterize the equality condition in Hu-Ye's conjecture for the eigenvalues of several hypergraph classes and present a new conjecture  regarding the eigenvalue multiplicities of  hypergraphs.
As an initial step toward this new conjecture, we give a sufficient condition for a point in the projective eigenvariety of a tensor to have local intersection multiplicity one.
\end{abstract}

\maketitle

\section{Introduction}
A \emph{tensor} (also called \emph{hypermatrix}) $\A=(a_{i_1\ldots i_k})$ of order $k$ and dimension $n$ over a field $\F$ is a multi-array of entries $ a_{i_1\ldots i_k} \in \F$ for all $i_j \in [n]:=\{1,\ldots,n\}$ and $j \in [k]$.
In particular, a square matrix corresponds to a tensor of order $2$.
For square matrices over the complex field $\C$, it is known that the algebraic multiplicity of any eigenvalue is always greater than or equal to its geometric multiplicity.
A natural generalization of this fact is to inquire about the relationship between the algebraic and geometric multiplicities of eigenvalues for general tensors.
To address this topic, we first recall some fundamental notions.
Given a vector $x \in \F^n$, the vector $\A x^{k-1} \in \F^n$ is defined component-wise as follows:
\begin{align*}
(\A x^{k-1})_i & =\sum_{i_{2},\ldots,i_{k}\in [n]}a_{ii_{2}\ldots i_{k}}x_{i_{2}}\cdots x_{i_k}, i \in [n].
\end{align*}
An $n$-dimensional tensor $\mathcal{I}=(i_{i_1i_2\ldots i_k})$ of order $k$ is called an \emph{identity tensor} if $i_{i_{1}i_2 \ldots i_{k}}=1$ for
   $i_{1}=i_2=\cdots=i_{k} \in [n]$, and $i_{i_{1}i_2 \ldots i_{k}}=0$ otherwise.

In 2005, Lim \cite{Lim} and Qi \cite{Qi} independently introduced tensor eigenvalues.

\begin{defi}[\cite{Lim,Qi}] Let $\A$ be a $k$-th order $n$-dimensional tensor over a field $\F$.
For some $\lambda \in \F$, if the polynomial system $(\lambda \mathcal{I}-\A)x^{k-1}=0$,
or equivalently, $\A x^{k-1}=\lambda x^{[k-1]}$, has a solution $x\in \F^{n}\backslash \{0\}$,
then $\lambda $ is called an \emph{eigenvalue} of $\A$ and $x$ is an \emph{eigenvector} of $\A$ associated with $\lambda$,
where $x^{[k-1]}:=(x_1^{k-1}, x_2^{k-1},\ldots,x_n^{k-1})$.
\end{defi}

The \emph{characteristic polynomial} $\varphi_\A(\la)$ of $\A$ is defined as the resultant of the polynomial system $(\la \I-\A)x^{k-1}$ (see \cite{Qi,CPZ09}).
It is known that $\la$ is an eigenvalue of $\A$ if and only if it is a root of $\varphi_\A(\la)$.
The \emph{algebraic multiplicity} of $\la$, denoted by $\am(\la)$, is the multiplicity of $\la$ as a root of $\varphi_\A(\la)$.
To avoid ambiguity, we occasionally write $\am(\la,\A)$ to emphasize the algebraic multiplicity of $\la$ as an eigenvalue of $\A$.
The \emph{spectrum} of $\A$, denoted by $\spec(\A)$, is the multiset of all the roots of $\varphi_\A(\la)$.
Unless stated otherwise, we work in the complex field $\F=\C$.
The \emph{spectral radius} of $\A$, denoted by $\rho(\A)$, is defined as the largest modulus of elements in $\spec(\A)$.
An eigenvalue of $\A$ is an \emph{H-eigenvalue} if it is associated with a real eigenvector.

Let $\la$ be an eigenvalue of a $k$-th order $n$-dimensional tensor $\A$ over $\C$.
The \emph{eigenvariety} of $\A$ associated with $\la$ is defined as the affine variety
$$\V_{\la}(\A)=\{x \in \C^n: \A x^{k-1}=\la x^{[k-1]}\},$$
which consists of all eigenvectors of $\A$ associated with $\la$ together with the zero vector.
The \emph{geometric multiplicity} of $\la$,  denoted by $\gm(\la)$, is defined as the dimension of $\V_\la(\A)$, or equivalently, the maximum dimension of the irreducible components of $\V_\la(\A)$.
The \emph{projective eigenvariety} of $\A$ associated with $\la$ is defined as the projective variety (\cite{FBH19})
$$\PV_\la(\A)=\{x \in \mathbb{P}^{n-1}: \A x^{k-1}=\la x^{[k-1]}\},$$
where $\mathbb{P}^{n-1}$ denotes the complex projective space of dimension $n-1$.

In 2016, Hu and Ye \cite{HuY16} proposed the following conjecture and established its validity in several cases.

\begin{conj}[\cite{HuY16} Hu-Ye Conjecture]
Suppose that a $k$-th order $n$-dimensional tensor $\A$ has an eigenvalue $\la$ with eigenvariety $\V_\la(\A)$ possessing $\kappa$ irreducible components $V_1,\ldots, V_\kappa$. Then
\begin{equation}\label{conj1} \am(\la) \ge \sum_{i=1}^\kappa \dim(V_i)(k-1)^{\dim(V_i)-1}.\end{equation}
In particular,
\begin{equation}\label{conj2} \am(\la) \ge \gm(\la)(k-1)^{\gm(\la)-1}.\end{equation}
\end{conj}

To verify the Hu-Ye Conjecture, one usually considers the adjacency tensor of a $k$-uniform hypergraph.
Cooper and Fickes \cite{CooperF22} fully described the irreducible components of the zero-eigenvariety of a $3$-uniform (loose) hyperpath $P_n^{(3)}$.
 Applying the result of Bao, Fan, Wang, and Zhu \cite{BaoFWZ20}, they calculated the algebraic multiplicity of the zero eigenvalue and confirmed the inequality \eqref{conj1} in Hu-Ye's conjecture for the zero eigenvalue of $P_n^{(3)}$.
Zheng \cite{Zheng24} confirmed the inequality \eqref{conj2} in Hu-Ye's conjecture for all nonzero eigenvalues of both the $k$-uniform hyperpath $P_n^{(k)}$ and the $k$-uniform hyperstar $S_n^{(k)}$, using the results of Chen and Bu \cite{ChenB21} and Bao, Fan, Wang, and Zhu \cite{BaoFWZ20} on the characteristic polynomials of hyperpaths and hyperstars.

Let $\A$ be a weakly irreducible nonnegative tensor (see Section \ref{PFsubsection} for the definition of weak irreducibility) with spectral radius $\rho$.
By the Perron-Frobenius theorem for nonnegative tensors, $\rho$ is an eigenvalue of $\A$.
Let $\la$ be any eigenvalue of $\A$ with modulus $\rho$.
We establish a relationship between the algebraic multiplicity $\am(\la)$ and the cardinality of the projective variety $\PV_\la(\A)$, which in turn implies that the Hu-Ye Conjecture holds for the eigenvalue $\la$ of $\A$.

\begin{thm}\label{amPvT}
Let $\A$ be a weakly irreducible nonnegative tensor of order $k$ and dimension $n$, and let $\la$ be any eigenvalue of $\A$ with modulus equal to the spectral radius of $\A$.
Then
$$ \am(\la) \ge |\PV_\la(\A)|,$$
and the Hu-Ye Conjecture holds for the eigenvalue $\la$ of $\A$.
\end{thm}

Recall that a real tensor $\B$ is called a \emph{$Z$-tensor} if all of its off-diagonal entries are non-positive; equivalently, it can be written as
$ \B=s \I - \A,$
where $s >0$ and $\A$ is nonnegative; see \cite{DingQW13, ZhangQZ14}.
The least H-eigenvalue of $\B$ is exactly $s-\rho(\A)$, which is also the eigenvalue of $\B$ with the least real part \cite{ZhangQZ14}.
The following result is an immediate consequence  of Theorem \ref{amPvT}.

\begin{coro}\label{amZ}
Let $\B$ be a weakly irreducible $Z$-tensor, and let $\la$ be the least H-eigenvalue of $\B$.
Then
$$ \am(\la) \ge |\PV_\la(\B)|,$$
and the Hu-Ye Conjecture holds for the eigenvalue $\la$ of $\B$.
\end{coro}

Let $H$ be a connected $k$-uniform hypergraph.
Its adjacency tensor $\A(H)$ and signless Laplacian tensor $\Q(H)$ are both nonnegative and weakly irreducible tensors, while its Laplacian tensor $\L(H)$ is a weakly irreducible $Z$-tensor.
These observations yield a broad family of tensors that satisfy Theorem \ref{amPvT} or Corollary \ref{amZ}.

Furthermore, for the adjacency tensors of certain classes of hypergraphs, we prove that the inequality in  Theorem \ref{amPvT} holds with equality. This, in turn, implies the equality case of inequality \eqref{conj1} in the Hu-Ye Conjecture.

\begin{thm}\label{amPv}
Let $H$ be one of the following classes of hypergraphs:
\begin{itemize}
\item[(1)] $k$-uniform hypertrees,
\item[(2)] $k$-th powers of connected simple graphs,
\item[(3)] complete $3$-uniform hypergraphs on at least $4$ vertices.
\end{itemize}
Then, for any eigenvalue $\la$ of the adjacency tensor $\A(H)$ with modulus equal to the spectral radius,
$$\am(\la)=|\PV_{\la}(\A(H))|.$$
Consequently, the equality case of inequality \eqref{conj1} in the Hu-Ye Conjecture holds for the eigenvalue $\la$.
\end{thm}

We conjecture that Theorem \ref{amPv} generalizes to  all connected uniform hypergraphs, leading to the following proposal.

\begin{conj}\label{conjH}
Let $H$ be a connected $k$-uniform hypergraph with spectral radius $\rho$.
Then $$\am(\rho)=|\PV_\rho(\A(H))|.$$
\end{conj}

It is known that zero is always an eigenvalue of the Laplacian tensor of a uniform hypergraph.
We have a similar result to Theorem \ref{amPv} for the zero Laplacian eigenvalue.

\begin{thm}\label{LamPv}
Let $T$ be a $k$-uniform hypertree, and let $\L(T)$ be the Laplacian tensor of $T$. Then
$$ \am(0, \L(T))= |\PV_0(\L(T))|.$$
Consequently, the equality case of the inequality \eqref{conj1} in the Hu-Ye Conjecture holds for the zero eigenvalue of $\L(T)$.
\end{thm}

Building on the proof of Theorem \ref{amPvT}, we initiate a discussion of Conjecture \ref{conjH}. As a foundational step toward resolving this conjecture, we establish a sufficient condition for an isolated point in a projective eigenvariety to have local intersection multiplicity $1$ (see Lemma \ref{simple2}). Using this result, we proceed to prove the simplicity of an arbitrary eigenvector $x \in \PV_\la(\A(H))$, where $|\la|=\rho(\A(H))$.

\begin{thm}\label{multirho}
Let $H$ be a connected $k$-uniform hypergraph.
Then for any eigenvalue $\la$ of $\A(H)$ with $|\la|=\rho(\A(H))$, each point $p \in \PV_{\la}(\A(H))$ has local intersection multiplicity $1$.
\end{thm}

The paper is organized as follows.
In Section 2 we will introduce some basic notions and preliminary results concerning nonnegative tensors and hypergraphs.
In Section 3,  we prove Theorem \ref{amPvT} and Corollary \ref{amZ} by applying the Poisson formula.
In Section 4, we prove Theorem \ref{amPv} and Theorem \ref{LamPv}, formulate Conjecture \ref{conjH}, and  discuss some related problems regarding the eigenvalues of Laplacian or signless Laplacian tensors.
Finally, in Section 5, we utilize the Jacobian matrix of polynomials to prove Theorem \ref{multirho} and several related results on the simplicity of the  eigenvectors associated with Laplacian or signless Laplacian tensors — findings that provide a foundational step towards resolving Conjecture \ref{conjH}.

\section{Preliminaries}
\subsection{Perron-Frobenius theorem for nonnegative tensors}\label{PFsubsection}
Let $\A=(a_{i_1\ldots i_k})$ be a tensor of order $k$ and dimension $n$.
Associate $\A$ with a directed graph $D(\A)$ on vertices $1,\ldots,n$ such that $(i, j)$ is an arc of $D(\A)$ if and only if there exists a nonzero entry $a_{i i_2\ldots i_k}$ such that $j \in \{i_2,\ldots,i_k\}$.
The tensor $\A$ is called \emph{weakly irreducible} if $D(\A)$ is strongly connected; otherwise, it is \emph{weakly reducible}.
We say $\A$ is \emph{nonnegative} if all its entries are nonnegative, and $\A$ is \emph{symmetric} if all its entries $a_{i_1\ldots i_k}$ are invariant under any permutation of its indices.
The \emph{support} (or \emph{zero-nonzero pattern} \cite{Shao13}) of $\A$, denoted by $\supp(\A)=(s_{i_1\ldots i_k})$, is the tensor of the same order and dimension as $\A$ defined by $s_{i_1\ldots i_k}=1$ if $a_{i_1\ldots i_k}\ne 0$, and $s_{i_1\ldots i_k}=0$ otherwise.
The tensor $\A$ is called \emph{combinatorially symmetric} if $\supp(\A)$ is symmetric \cite{FBH19}.
Surely, combinatorially symmetric tensors are generalizations of symmetric tensors.
For a combinatorially symmetric tensor $\A$ of order $k$ and dimension $n$, define
$$ E(\A)=\{(i_1,\ldots,i_k) \in [n]^k: a_{i_1\ldots i_k} \ne 0, i_1\le \cdots \le i_k\}.$$
For each $e=(i_1,\ldots,i_k)\in E(\A)$ and $j \in [n]$, let
$$ b_{e,j}=|\{t: i_t=j, t \in [k]\}|.$$
The matrix $B_\A=(b_{e,j})$ is called the \emph{incidence matrix} of $\A$ \cite{FBH19}.

Chang et al. \cite{CPZ08}, Yang and Yang \cite{YangY10,YangY11-1,YangY11-2}, and Friedland et al. \cite{FriGH} generalized the Perron-Frobenius theorem from nonnegative matrices to nonnegative tensors.
In the following, we list part of the results of the theorem,
where an eigenvalue is called an \emph{$H^{++}$-eigenvalue} if it is associated with a positive eigenvector.

\begin{thm}[\cite{FriGH,YangY11-2}]\label{PFthm}
Let $\A$ be a weakly irreducible nonnegative tensor of order $k$ and dimension $n$ with spectral radius $\rho(\A)$. Then the following results hold.

\begin{itemize}
\item[(1)] $\rho(\A)$ is a unique  $H^{++}$-eigenvalue, with a unique positive eigenvector up to a positive scalar, called the \emph{Perron vector} of $\A$.

\item[(2)]  If $x$ is an eigenvector of $\A$ associated with an eigenvalue with modulus $\rho(\A)$, then $|x|$ is the Perron vector of $\A$.

\item[(3)] If $\A$ has $m$ distinct eigenvalues with modulus equal to $\rho(\A)$, then these eigenvalues are $\rho(\A)e^{\i 2 \pi j/m}, j=0,1,\ldots,m-1$. Furthermore, there exists a diagonal matrix $D$ with unit-modulus diagonal entries such that $\A=e^{-\i 2\pi /m}D^{-(k-1)}\A D$, and $\spec(\A)=e^{\i 2\pi /m}\spec(\A)$.
\end{itemize}

\end{thm}

In Theorem \ref{PFthm}, $\i=\sqrt{-1}$.
For two $n\times n$ diagonal matrices $P$ and $Q$, the product $P\A Q$ is defined by Shao \cite{Shao13} as
$$(P\A Q)_{i_1,\ldots,i_k}=P_{i_1 i_1} a_{i_1,\ldots,i_k} Q_{i_2 i_2} \cdots Q_{i_k i_k}.$$
If $D$ is an invertible diagonal matrix, then $D^{-(k-1)}\A D$ is said to be \emph{diagonally similar to} $\A$, and has the same spectrum as $\A$ \cite{Shao13}.

\subsection{Cyclic index}
The value $m$ in Theorem \ref{PFthm}(3) is called the \emph{cyclic index} of $\A$ (\cite{CPZ11}), and is related to the spectral symmetry of $\A$.

\begin{defi}[\cite{FHB19}]\label{sym}
Let $\A$  be a general tensor, and let $\ell$ be a positive integer. The tensor $\A$
 is called \emph{spectral $\ell$-symmetric} if
 \begin{equation}\label{ell}
 \spec(\A)=e^{\i 2\pi /\ell}\spec(\A).
 \end{equation}
 The \emph{cyclic index} of $\A$, denoted by $c(\A)$, is the maximum positive integer $\ell$ for which \eqref{ell} holds.
\end{defi}

If $\A$ is nonnegative and weakly irreducible, the cyclic index $c(\A)$ in Definition \ref{sym} is consistent with that introduced in \cite{CPZ09}; see \cite{FHB19} for further details.
Using the following result,
Fan et al.  \cite{FHB19} proved that
 $ c(\A)=\gcd \{d: \Tr_d(\A) \ne 0\},$
where $\Tr_d(\A)$ is the generalized $d$-th order trace of $\A$ (see \cite{CooperD12,MS10,ShaoQH15}).

\begin{thm}[\cite{FHB19}]\label{symequ}
Let $\A$ be a tensor of order $k$ and dimension $n$.
Let $\varphi_\A(\la)=\sum_{i=0}^N a_i \la^{N-i}$ be the characteristic polynomial of $\A$, where $N=n(k-1)^{n-1}$. Then the following conditions are equivalent.

\begin{itemize}
\item[(1)] $\A$ is spectral $\ell$-symmetric.

\item[(2)] If $\ell \nmid d$, then $a_d =0$; that is, $\varphi_\A(\la)=\la^t f(\la^\ell)$ for some nonnegative integer $t$ and some polynomial $f$.

\item[(3)] If $\ell \nmid d$, then $\Tr_d(\A) =0$.

\end{itemize}
\end{thm}

\begin{coro}\label{multi-coro}
Let $\A$ be a weakly irreducible nonnegative tensor.
Then, for any eigenvalue $\la$ of $\A$ with modulus $\rho(\A)$,
$\am(\la)=\am(\rho(\A)).$
\end{coro}

\begin{proof}
Let $\la$ be an eigenvalue of $\A$ with modulus $\rho(\A)$.
By Theorem \ref{PFthm}(3), $\la=\rho(\A)e^{\i 2 \pi j /m}$ for some $j \in \{0,1,\ldots,m-1\}$, where $m=c(\A)$.
Furthermore, $\A$ is spectral $m$-symmetric.
By Theorem \ref{symequ}(2),
the characteristic polynomial of $\A$ satisfies $\varphi_\A(\la)=\la^t f(\la^m)$  for some nonnegative integer $t$ and some polynomial $f$.
Therefore, $\la$ and $\rho(\A)$ have the same algebraic multiplicity.
\end{proof}

\subsection{Stabilizing index}
Next, we introduce the stabilizing index of a tensor, which is related to the projective eigenvariety associated with the spectral radius.
Let $\A$ be a weakly irreducible nonnegative tensor of order $k$ and dimension $n$, which is also spectral $\ell$-symmetric.
 Let $$\la_j=\rho(\A)e^{\i 2 \pi j /\ell}, j=0,1,\ldots,\ell-1$$
 be the eigenvalues of $\A$, where $\la_0=\rho(\A)$.
Define
$$ \mathbb{D}(\A)=\bigcup_{j=0}^{\ell-1} \mathbb{D}_j(\A),~ \mathbb{D}_j(\A)=\{D: \A=e^{\i 2 \pi j/\ell} D^{-(k-1)}\A D, d_{11}=1\},j=0,1,\ldots,\ell-1,$$
where $D$ is an invertible $n \times n$ diagonal matrix  in the above definition.
Note that $\mathbb{D}_0(\A) \ne \emptyset$ as it contains the identity matrix.
It is straightforward to verify that $\mathbb{D}_0(\A)$ forms an abelian group under standard matrix multiplication.

\begin{defi}[\cite{FBH19}]
The \emph{stabilizing index} of a tensor $\A$, denoted by $s(\A)$, is defined to be the cardinality of the group $\mathbb{D}_0(\A)$.
\end{defi}

Let
$\PV(\A):=\cup_{j=0}^{\ell-1} \PV_{\la_j}(\A)$.
By Theorem \ref{PFthm}(2), each eigenvector $x \in \PV(\A)$ has no zero entries, so we can assume that each point $x \in \PV(\A)$ holds $x_1=1$ by normalization.
Write $x=D_{x} \p$, where $\p:=|x|$ is the unique Perron vector of $\A$ also by Theorem \ref{PFthm}(2), and
$$ D_{x}=\diag(x_1/|x_1|, \ldots, x_n/|x_n|).$$
We equip $\PV(\A)$ with a binary operation $\circ$ defined by
$ x \circ x' := D_{x} D_{x'} \p,$
and define a map:
\begin{equation}\label{Psi} \Psi: \PV(\A) \to \mathbb{D}(\A), x \mapsto D_{x}.
\end{equation}

\begin{lem}[\cite{FBH19}]\label{group}
Let $\A$ be a weakly irreducible nonnegative tensor of order $k$ and dimension $n$ that is spectral $\ell$-symmetric. The following results hold.

\begin{itemize}

\item[(1)] $(\mathbb{D}(\A), \cdot)$ is an abelian group under standard matrix multiplication, where $\mathbb{D}_0(\A)$ is a subgroup of
$\mathbb{D}(\A)$, and $\mathbb{D}_j(\A)$ is a coset of $\mathbb{D}_0(\A)$ for $j \in [\ell-1]$.

\item[(2)] $(\PV(\A),\circ)$ is an abelian group, where $\PV_{\rho(\A)}(\A)$ is a subgroup of $\PV(\A)$, and $\PV_{\la_j}(\A)$ is a coset of $\PV_{\rho(\A)}(\A)$ for each $j \in [\ell-1]$. Moreover, $\Psi$ is a group isomorphism from $\PV(\A)$ to $\mathbb{D}(\A)$ which sends $\PV_{\la_j}(\A)$ to $\mathbb{D}_j(\A)$ for each $j=0,1,\ldots,\ell-1$.

\item[(3)] If $\A$ is additionally combinatorially symmetric, then $\ell \mid k$ and $D^k=I$ for every $D \in \mathbb{D}(\A)$.

\end{itemize}
\end{lem}

Fan et al. \cite{FanHB22} proved that $\PV_{\rho(\A)}(\A)$ is finite for any weakly irreducible nonnegative tensor $\A$. Combining this with Lemma \ref{group}, we deduce the following result.

\begin{coro}[\cite{FanHB22}]\label{Pv-spectral}
Let $\A$ be a weakly irreducible nonnegative tensor of order $k$ and dimension $n$ that is spectral $\ell$-symmetric.
For each $j \in [\ell-1]$, let $\la_j=\rho(\A)e^{\i 2 \pi j /\ell}$.
Then $|\PV_{\la_j}(\A)|=|\PV_{\rho(\A)}(\A)|=s(\A)$.
\end{coro}

Moreover, assume that $\A$ is combinatorially symmetric.
By Lemma \ref{group}(3), for any $y \in \PV(\A)$, $y^{\circ k}=(D_y)^k \p=\p$, so $\PV(\A)$ and $\mathbb{D}(\A)$ both admit $\Z_k$-modules.
By Lemma \ref{group}(2), $\PV(\A)$ is $\Z_k$-isomorphic to $\mathbb{D}(\A)$, which sends $\PV_{\la_j}(\A)$ to $\mathbb{D}_j(\A)$ for each $j=0,1,\ldots,\ell-1$.
For any $D \in \mathbb{D}(\A)$, as $D^k=I$, we can write $$D=\diag(e^{\i 2 \pi \phi_1(D)/k},\ldots, e^{\i 2 \pi \phi_n(D)/k}),$$
 where $\phi_i(D) \in \Z_k$ for $i \in [n]$, and $\phi_1(D)=0$ as $D_{11}=1$ by definition.
 Thus, we define the map
\begin{equation}\label{Phi}\Phi: \mathbb{D}(\A) \to \Z_k^n, D \mapsto (\phi_1(D), \ldots, \phi_n(D)).
\end{equation}
Let $\mathbf{1}$ denote the all-one vector of dimension $n$.
By Lemma \ref{group}(3), $\ell$ divides $k$, so $kj/\ell$ is an integer.
For each $j=0,1,\ldots,\ell-1$, define
\begin{equation}\label{PS}
\S_j(\A)=\{x \in \Z_k^n: B_\A x= \frac{kj}{\ell} \mathbf{1} \mbox{~over~} \Z_k, x_1=0\}, \end{equation}
and set $\S(\A)=\cup_{j=0}^{\ell-1} \S_j(\A)$.
Clearly, $\S(\A)$ is a $\Z_k$-module that contains $\S_0(\A)$ as a $\Z_k$-submodule, and for each $j \in [\ell-1]$, $\S_j(\A)$ is a coset of $\S_0(\A)$, i.e., $\S_j(\A)=\S_0(\A)+y^{(j)}$ for some $y^{(j)} \in \S_j(\A)$.

\begin{lem}[\cite{FBH19}]\label{iso}
Let $\A$ be a weakly irreducible nonnegative combinatorially symmetric tensor that is spectral $\ell$-symmetric. Then $\Phi$ is a $\Z_k$-isomorphism from $\mathbb{D}(\A)$ to $\S(\A)$ which sends $\mathbb{D}_j(\A)$ to $\S_j(\A)$ for each $j=0, 1, \ldots, \ell-1$.
\end{lem}

By Lemmas \ref{group} and \ref{iso}, $\Phi \Psi: \PV(\A) \to \S(\A)$ is a $\Z_k$-isomorphism sending $ \PV_{\la_j}(\A)$ to $\S_j(\A)$ for $j=0,1,\ldots, \ell-1$.
So, for each $x \in \PV_{\la_j}(\A)$, writing
$$x=(|x_1| e^{\i 2 \pi \phi_1(x)/k}, \ldots, |x_n| e^{\i 2 \pi \phi_n(x)/k}),$$
 where $\phi_i(x) \in \Z_k$ for $i \in [n]$.
We have
$(\phi_1(x), \ldots, \phi_n(x)) \in \S_j(\A)$, namely,
for each nonzero entry $a_{i_1\ldots i_k}$ of $\A$,
\begin{equation}\label{angle}
\phi_{i_1}(x)+\cdots+\phi_{i_k}(x) \equiv \frac{kj}{\ell} \mod k.
\end{equation}
We refer to \eqref{angle} as the \emph{uniform angle property}.

The following theorem characterizes the $\Z_k$-module structure of $\S_0(\A)$, thus determining the stabilizing index $s(\A)$ of $\A$.

\begin{thm}[\cite{FBH19}]\label{PV-Zk}
Let $\A$ be a weakly irreducible nonnegative combinatorially symmetric tensor of order $k$ and dimension $n$.
Suppose that Smith normal form of $B_\A$ over $\Z_k$ has invariant divisors $d_1,\ldots,d_r$. Then $1 \le r \le n-1$, and as $\Z_k$-modules,
$$ \PV_{\rho(\A)}(\A)\cong \S_0(\A) \cong \oplus_{i,d_i \ne 1} \Z_{d_i} \oplus \underbrace{\Z_k \oplus \cdots \oplus \Z_k}_{n-r-1 \mbox{~copies}}.$$
Consequently, the stabilizing index of $\A$ satisfies $s(\A)=|\PV_{\rho(\A)}(\A)|= k^{n-r-1}\prod_{i=1}^r d_i$.
\end{thm}

\subsection{Hypergraphs and tensor representations}
A \emph{hypergraph} $H=(V(H),E(H))$ consists of a vertex set $V(H)$ and an edge set $E(H)$ whose elements are the subsets of $V(H)$.
The \emph{degree} of a vertex $v$ of $H$, denoted by $d_v$, is the number of edges of $H$ containing $v$.
If, for each edge $e \in E(H)$, $|e|=k$, then $H$ is called \emph{$k$-uniform}.

Let $H$ be a $k$-uniform hypergraph on vertices $v_1,\ldots,v_n$.
Cooper and Dutle \cite{CooperD12} introduced the \emph{adjacency tensor} $\mathcal{A}(H)=(a_{i_1 i_2 \ldots i_k})$ of $H$, a $k$-th order $n$-dimensional tensor defined by
        \[  a_{i_{1}i_{2}\dots i_{k}}=\begin{cases}
    	\frac{1}{(k-1 )!},  & \text{ if } \{ v_{i_{1} },v_{i_{2} },\dots v_{i_{k} }    \}\in E(H ),   \\
    	~~0,  & \text{ otherwise }.
    \end{cases}      \]
The \emph{Laplacian tensor} $\L(H)$ and the \emph{signless Laplacian tensor} $\Q(H)$
are defined as $\L(H)=\D(H)-\A(H)$ and $\Q(H)=\D(H)+\A(H)$, respectively \cite{Qi14}. Here, $\D(H)$ is a diagonal tensor of order $k$ and dimension $n$ with the $(v,v,\ldots,v)$-entry equal to $d_v$ (the degree of the vertex $v$) for each $v \in V(H)$.
It is known that $\A(H),\Q(H)$ are nonnegative,  $\A(H),\L(H),\Q(H)$ are all symmetric, and each of these tensors is weakly irreducible if and only if $H$ is connected \cite{YangY11-2}.

The \emph{incidence matrix} of $H$, denoted by $B_H=(b_{e,v})$, coincides with the incidence matrix of $\A(H)$, that is, $b_{e,v}=1$ if $v \in e$, and $b_{e,v}=0$ otherwise.
When viewing $B_H$ as a matrix over $\Z_k$, the incidence matrices of $\L(H),\Q(H)$ also coincide with $B_H$.

Let $H$ be a connected $k$-uniform hypergraph on $n$ vertices.
Define
\begin{equation} \label{S0-def} \S_0(H)=\{x \in \Z_k^n: B_H x =\mathbf{0}, x_1=0\},
\end{equation}
where $\mathbf{0}$ denotes the zero vector in $\Z_k^n$.
By Theorem \ref{PV-Zk}, we have
\begin{equation}\label{Sadj} |\PV_{\rho(\A(H))}(\A(H))|=|\PV_{\rho(\Q(H))}(\Q(H))|=|\S_0(H)|.
\end{equation}

It is known that $0$ is always an eigenvalue of $\L(H)$ associated with the all-one vector $\mathbf{1}$.
Hu and Qi \cite{HuQ14} established the following properties for the eigenvectors $x$ of $\L(H)$ associated with the zero eigenvalue:
(1) $x$ contains no zero entries, and (2) if normalizing $x$ such that $x_1=1$, then $x$ has the form
\begin{equation}\label{Lzero} x=(e^{\i 2 \pi \phi_1(x)/k}, \ldots, e^{\i 2 \pi \phi_n(x)/k}), \end{equation}
 where $\phi_i(x) \in \Z_k$ and $\phi_1(x)=0$.
 Furthermore, $x$ satisfies the uniform angle property, namely, for each edge $e \in E(H)$,
\begin{equation}\label{unimod} \sum_{i \in e} \phi_i(x)= 0 \mod k.\end{equation}
Thus, after normalization, each point of
\(\PV_0(\L(H))\) has the form \eqref{Lzero}.
 Fan et al. \cite{FanW19} introduced a Hadamard product in $\PV_0(\L(H))$, and proved that $\PV_0(\L(H))$ is a $\Z_k$-module isomorphic to $\S_0(H)$.
Consequently,
\begin{equation}\label{Slap}|\PV_0(\L(H))|=|\S_0(H)|.
\end{equation}

We conclude with a discussion of the zero eigenvalue of $\Q(H)$.
It was proved that zero is an eigenvalue of $\Q(H)$ if and only if $k$ is even and $H$ is odd-colorable \cite{Niki, HuQ14,FanW19}.
For any eigenvector $x$ of $\Q(H)$ associated with the zero eigenvalue,
normalizing $x$ to satisfy $x_1=1$ implies that $x$ also has the form \eqref{Lzero} and satisfies the uniform angle property (\cite{HuQ14}): for each $e \in E(H)$,
\begin{equation}\label{Qunimod} \sum_{i \in e} \phi_i(x)= \frac{k}{2} \mod k.\end{equation}
Fan et al. \cite{FanW19} further showed that if zero is an eigenvalue of $\Q(H)$, then
\begin{equation}\label{Sslap} |\PV_0(\Q(H))|=|\PV_0(\L(H))|=|\S_0(H)|.
\end{equation}

\section{Hu-Ye Conjecture on spectral radius of tensors}
Let $\K$ be an algebraically closed field, and let $F_1, \ldots, F_n \in \K[x_1,\ldots,x_n]$ be homogeneous polynomials of degrees $d_1,\ldots, d_n$, respectively.
Define
\begin{equation}\label{GenH0} \bar{F}_i(x_1,\ldots,x_{n-1}):=F_i(x_1,\ldots,x_{n-1},0), i \in [n],
\end{equation}
\begin{equation}\label{GenH1}  f_i(x_1,\ldots,x_{n-1}):=F_i(x_1,\ldots,x_{n-1},1), i \in [n].
\end{equation}

\begin{thm}[Poisson formula; cf. \cite{Cox}, Chapter 3]\label{poisson}
If $\Res(\bar{F}_1,\ldots,\bar{F}_{n-1})\ne 0$, then the quotient ring
$A = \K[x_1,\ldots,x_{n-1}]/\langle f_1,\ldots,f_{n-1}\rangle$ is a finite-dimensional $\K$-vector space of dimension $d_1\cdots d_{n-1}$. Moreover,
$$ \Res(F_1,\ldots,F_n)=\Res(\bar{F}_1,\ldots,\bar{F}_{n-1})^{d_n} \det(m_{f_n}: A \to A),$$
where $m_{f_n}: A \to A$ is the $\K$-linear map given by multiplication by $f_n$.
\end{thm}

It is known that for each $p \in \V(\langle f_1,\ldots,f_{n-1}\rangle)$, the value $f_n(p)$ is an eigenvalue of $m_{f_n}$ with multiplicity $m(p)$ (the local intersection multiplicity defined in Section \ref{MultiZero}), and
\begin{equation}\label{deter} \det m_{f_n} = \prod_{p \in \V(\langle f_1,\ldots,f_{n-1}\rangle)} f_n(p)^{m(p)},\end{equation}
a result that can be found in \cite{Cox}, Chapter 4.

We first recall two standard facts used in the specialization argument.  Let $R$ be a discrete valuation ring with uniformizer $t$, fraction field $K$, and algebraically closed residue field $\mathbf{k}=R/(t)$.  Every $0\ne a\in K$ can be written uniquely as
\[
 a=t^{v_t(a)}u,\quad u\in R^\times,
\]
which defines the normalized $t$-adic valuation $v_t:K^\times\to\mathbb Z$; we set $v_t(0)=+\infty$.

We also recall the local decomposition of a finite-dimensional commutative algebra over $\mathbf{k}$.  Given polynomials $f_1^{(0)},\ldots,f_l^{(0)}\in\mathbf{k}[x_1,\ldots,x_l]$, put
\[
 I_0:=\langle f_1^{(0)},\ldots,f_l^{(0)}\rangle,
 \quad
 Q_0:=\mathbf{k}[x_1,\ldots,x_l]/I_0,
\]
and suppose that $Q_0$ is finite-dimensional.  Let $q=(q_1,\ldots,q_l)$ be a point of $\V(I_0)$, and set
\[
 M_q:=\langle x_1-q_1,\ldots,x_l-q_l\rangle
 \subseteq\mathbf{k}[x_1,\ldots,x_l].
\]
Since $q\in\V(I_0)$, we have $I_0\subseteq M_q$.  The evaluation map and its kernel are
\[
 \operatorname{ev}_q:Q_0\longrightarrow\mathbf{k},\quad [g]\longmapsto g(q),
\]
and
\[
 \mathfrak m_q:=\ker(\operatorname{ev}_q)
 =M_q/I_0
 =\{[g]\in Q_0:g(q)=0\}
 =\langle \bar x_1-q_1,\ldots,\bar x_l-q_l\rangle_{Q_0},
\]
where $\bar x_i$ is the residue class of $x_i$ in $Q_0$.  Thus $\mathfrak m_q$ is the maximal ideal corresponding to $q$.  The localization of $Q_0$ at $\mathfrak m_q$ is
\[
 A_q:=(Q_0)_{\mathfrak m_q}
 =\left\{
 \frac{[g]}{[h]}:
 g,h\in\mathbf{k}[x_1,\ldots,x_l],\ h(q)\ne0
 \right\}.
\]
Equivalently, since localization commutes with taking quotients,
\[
 A_q\cong
 \frac{\mathbf{k}[x_1,\ldots,x_l]_{M_q}}
 {I_0\mathbf{k}[x_1,\ldots,x_l]_{M_q}}.
\]
Thus $A_q$ is the local algebra of $Q_0$ at $q$.  Since $Q_0$ is finite-dimensional over $\mathbf{k}$, it is Artinian and admits the decomposition
\begin{equation}\label{Artinian-decomposition}
 Q_0\cong
 \prod_{q\in\V(I_0)}A_q,
\end{equation}
where every $A_q$ is an Artinian local $\mathbf{k}$-algebra and its maximal ideal $\mathfrak m_qA_q$ is nilpotent; see \cite[Chapters~3 and~8, especially Theorem~8.7]{AtiyahMacdonald}.  For $b\in Q_0$, let $T_b:Q_0\to Q_0$ denote multiplication by $b$.  For each $q$, let
\[
 \iota_q:Q_0\longrightarrow A_q,\quad [g]\longmapsto [g]/1,
\]
be the canonical localization map, write $b_q:=\iota_q(b)$, and define
\[
 T_{b,q}:A_q\longrightarrow A_q,\quad u\longmapsto b_qu.
\]
Equivalently, $T_{b,q}$ is the localization of $T_b$ at $\mathfrak m_q$.  Under the decomposition in \eqref{Artinian-decomposition}, $T_b$ decomposes as the product of the maps $T_{b,q}$, and hence
\begin{equation}\label{Artinian-kernel}
 \dim_{\mathbf{k}}\ker T_b
 =\sum_{q\in\V(I_0)}
   \dim_{\mathbf{k}}\ker T_{b,q}.
\end{equation}

The following lemma is the form of the Poisson formula needed below.  It works over the local parameter ring and therefore does not require choosing algebraic branches of the points of a generic fiber.

\begin{lem}\label{specialization-lemma}
Let $R,t,K,\mathbf{k}$, and $v_t$ be as above.  Let
$F_1,\ldots,F_{l+1}\in R[x_1,\ldots,x_l,z]$ be homogeneous polynomials of positive degrees $d_1,\ldots,d_{l+1}$, respectively, and set
\[
 \bar F_i(x_1,\ldots,x_l)=F_i(x_1,\ldots,x_l,0),\quad
 f_i(x_1,\ldots,x_l)=F_i(x_1,\ldots,x_l,1).
\]
Suppose that $r:=\Res(\bar F_1,\ldots,\bar F_l)$ is a unit of $R$ and that $\Res(F_1,\ldots,F_{l+1})\ne 0$.  Then
\[
 Q=R[x_1,\ldots,x_l]/\langle f_1,\ldots,f_l\rangle
\]
is a finite free $R$-module of rank $d_1\cdots d_l$.  Let $M$ be multiplication by $f_{l+1}$ on $Q$, and let $M_0$ be its reduction on
\[
 Q_0=Q/tQ=\mathbf{k}[x_1,\ldots,x_l]/\langle f_1^{(0)},\ldots,f_l^{(0)}\rangle,
\]
where $f_i^{(0)}$ denotes reduction modulo $t$.  Then
\begin{equation}\label{specialization-bound}
 v_t\!\left(\Res(F_1,\ldots,F_{l+1})\right)
 =v_t(\det M)
 \ge \dim_{\mathbf{k}}\ker M_0
 \ge |Z_0|,
\end{equation}
where
\[
 Z_0=\{q\in\mathbf{k}^l:f_1^{(0)}(q)=\cdots=f_{l+1}^{(0)}(q)=0\}.
\]
\end{lem}

\begin{proof}
Since $r$ is a unit, its images in $K$ and in $\mathbf{k}$ are nonzero.  The defining property of the resultant therefore shows that the equations $\bar F_1=\cdots=\bar F_l=0$ have no common projective zero in either the generic or the closed fiber.  Equivalently, the homogenized equations $F_1=\cdots=F_l=0$ have no point on the hyperplane $z=0$ in either fiber.  Let
\[
 X:=V_+\!\left(\langle F_1,\ldots,F_l\rangle\right)
 =\operatorname{Proj}\frac{R[x_1,\ldots,x_l,z]}
 {\langle F_1,\ldots,F_l\rangle}
 \subseteq\mathbb P_R^l.
\]
Then $X\cap V(z)=\emptyset$, so the standard affine chart $D_+(z)$ contains all of $X$, and
\[
 X=X\cap D_+(z)
 \cong\operatorname{Spec}Q
 =V\!\left(\langle f_1,\ldots,f_l\rangle\right)
 \subseteq\mathbb A_R^l.
\]
The scheme $X$ is projective over $R$, as it is closed in a projective space over $R$, and the displayed identification shows that it is affine over $R$.  An affine and projective morphism is finite; consequently, $Q$ is a finite $R$-module.

Put $d=d_1\cdots d_l$, and write $Q_K:=Q\otimes_R K$ for the scalar extension of $Q$ from $R$ to its fraction field $K$; equivalently, $Q_K$ is the coordinate ring of the generic fiber of $\operatorname{Spec}Q\to\operatorname{Spec}R$.  On the generic fiber, the resultant of the leading homogeneous parts is nonzero, so the first assertion of Theorem \ref{poisson} gives
\[
 \dim_K Q_K=d.
\]
The same argument on the closed fiber, where the reduction of $r$ remains nonzero, gives
\[
 \dim_{\mathbf{k}}Q_0=d.
\]
By the structure theorem for finite modules over a discrete valuation ring, there are integers $b_1,\ldots,b_c\ge1$ such that
\[
 Q\cong R^d\oplus\bigoplus_{j=1}^{c}R/(t^{b_j}).
\]
After reduction modulo $t$, the right-hand side has $\mathbf{k}$-dimension $d+c$.  The preceding equality therefore forces $c=0$, and hence $Q$ is free of rank $d$.

Choose an $R$-basis of $Q$, so that $M$ is represented by a $d\times d$ matrix over $R$.  Applying Theorem \ref{poisson} after extending $K$ to its algebraic closure gives
\[
 \Res(F_1,\ldots,F_{l+1})=r^{d_{l+1}}\det M
\]
in $K$; the matrix of multiplication after this scalar extension is still the base change of $M$.  Since $r$ is a unit, this identity yields the first equality in \eqref{specialization-bound}.

The hypothesis $\Res(F_1,\ldots,F_{l+1})\ne0$ and the displayed Poisson identity show that $\det M\ne0$.  Thus the Smith normal form over $R$ has the shape
\[
 UMV=\operatorname{diag}(t^{r_1},\ldots,t^{r_s},1,\ldots,1),
 \quad U,V\in\operatorname{GL}_d(R),\quad r_i\ge1.
\]
Upon reducing modulo $t$, the matrices $U$ and $V$ remain invertible, and the displayed diagonal matrix has a kernel of dimension exactly $s$.  Consequently,
\[
 v_t(\det M)=\sum_{i=1}^s r_i\ge s=\dim_{\mathbf{k}}\ker M_0.
\]

For $q\in Z_0$, let $[f_{l+1}^{(0)}]$ denote the residue class of $f_{l+1}^{(0)}$ in $Q_0$, and put
\[
 a_q:=\iota_q([f_{l+1}^{(0)}])=[f_{l+1}^{(0)}]/1\in A_q.
\]
Since $f_{l+1}^{(0)}(q)=0$, we have $a_q\in\mathfrak m_qA_q$.  The nilpotence of this maximal ideal implies that multiplication by $a_q$ is a nilpotent endomorphism of the nonzero vector space $A_q$, and hence has a nonzero kernel.  Applying \eqref{Artinian-kernel} with $b=[f_{l+1}^{(0)}]$ gives
\[
 \dim_{\mathbf{k}}\ker M_0
 \ge\sum_{q\in Z_0}1
 =|Z_0|.
\]
\end{proof}

Let $\A=(a_{i_1\ldots i_k})$ be a  tensor over $\C$ of order $k$ and dimension $n$.
For each $i \in [n]$, define the polynomial
\begin{equation}\label{eigV} F_i(x_1,\ldots,x_n)=\la x_i^{k-1}-\sum_{i_2,\ldots,i_k}a_{i i_2 \ldots i_k}x_{i_2}\cdots x_{i_k}.
\end{equation}
Then the characteristic polynomial of $\A$ is given by the resultant
$$ \varphi_{\A}(\la)=\Res (F_1,\ldots,F_n).$$

\begin{proof}[\bf Proof of Theorem \ref{amPvT}]
Let $\la_0$ be an eigenvalue of $\A$ with $|\la_0|=\rho:=\rho(\A)$.  By Corollary \ref{multi-coro}, $\am(\la_0)=\am(\rho)$, and by Corollary \ref{Pv-spectral}, $|\PV_{\la_0}(\A)|=|\PV_{\rho}(\A)|$.  It therefore suffices first to prove the desired inequality for $\rho$.

The case $n=1$ is immediate, so assume $n\ge2$.

In the remainder of the proof, $\la$ denotes the indeterminate in the characteristic polynomial.  For each $i\in[n]$, let $F_i$ be as in \eqref{eigV}, and let $\bar F_i$ and $f_i$ be given by \eqref{GenH0} and \eqref{GenH1}.  Let $\B$ be the principal subtensor of $\A$ indexed by $[n-1]$.  Then
\[
 h(\la):=\Res(\bar F_1,\ldots,\bar F_{n-1})=\varphi_{\B}(\la).
\]
Since $\A$ is nonnegative and weakly irreducible, \cite[Corollary 3.4]{KhanF} gives $\rho(\B)<\rho(\A)=\rho$.  Hence $h(\rho)\ne0$.

Set $t=\la-\rho$ and $R=\C[\la]_{(t)}$.  The polynomial $h$ is a unit of the discrete valuation ring $R$.  
Moreover, $\Res(F_1,\ldots,F_n)=\varphi_{\A}(\la)$ is a
nonzero polynomial, so all hypotheses of Lemma \ref{specialization-lemma}
are satisfied.
Apply Lemma \ref{specialization-lemma}, with $l=n-1$, $z=x_n$, and all degrees equal to $k-1$, to the polynomials $F_1,\ldots,F_n$.  If
\[
 Q=R[x_1,\ldots,x_{n-1}]/\langle f_1,\ldots,f_{n-1}\rangle
\]
and $M$ is multiplication by $f_n$ on $Q$, then
\begin{equation}\label{rho-kernel-bound}
 \am(\rho)=v_t(\varphi_\A)=v_t(\det M)
 \ge \dim_{\C}\ker M_0\ge |Z_\rho|,
\end{equation}
where
\[
 Z_\rho=\{q\in\C^{n-1}:f_1(q)|_{\la=\rho}=\cdots=f_n(q)|_{\la=\rho}=0\}.
\]

Every point of $\PV_\rho(\A)$ has no zero coordinates by Theorem \ref{PFthm}(2).  It therefore has a unique representative whose last coordinate is $1$, and the map
\[
 Z_\rho\longrightarrow\PV_\rho(\A),\quad
 (q_1,\ldots,q_{n-1})\longmapsto[q_1:\cdots:q_{n-1}:1],
\]
is a bijection.  It follows from \eqref{rho-kernel-bound} that
\[
 \am(\rho)\ge |Z_\rho|=|\PV_\rho(\A)|.
\]

Returning to $\la_0$, the two equalities at the beginning of the proof give
\[
 \am(\la_0)\ge |\PV_{\la_0}(\A)|.
\]
The projective eigenvariety $\PV_{\la_0}(\A)$ is finite by Corollary \ref{Pv-spectral}.  Write it as $\{[x^{(1)}],\ldots,[x^{(\kappa)}]\}$, where the $x^{(i)}$ are nonzero representatives and $\kappa=|\PV_{\la_0}(\A)|$.  Since the eigenvector equations are homogeneous, the irreducible components of the affine eigenvariety $\V_{\la_0}(\A)$ are precisely the lines $\C x^{(1)},\ldots,\C x^{(\kappa)}$, each of dimension $1$.  Therefore,
\[
 \am(\la_0)\ge |\PV_{\la_0}(\A)|
 =\sum_{i=1}^{\kappa}\dim(\C x^{(i)})(k-1)^{\dim(\C x^{(i)})-1},
\]
which verifies \eqref{conj1}; in particular, $\gm(\la_0)=1$ and \eqref{conj2} also holds.
\end{proof}

\begin{proof}[\bf Proof of Corollary \ref{amZ}]
Let $\B=s \I - \A$,
where $s >0$ is a positive constant, and $\A$ is a weakly irreducible nonnegative tensor.
We first note that the algebraic multiplicity $\am(\la)$ of the eigenvalue $\la$ of $\B$ coincides with the algebraic multiplicity of $\rho(\A)$ as an eigenvalue of $\A$.
Observe that $\PV_\la(\B)=\PV_{\rho(\A)}(\A)$ as a vector $x$ is an eigenvector of $\B$ associated with $\la$ if and only if $x$ is an
 eigenvector of $\A$ associated with $\rho(\A)$ (cf. \cite[Lemma 3.1]{FanWB19}).
The desired result then follows from Theorem \ref{amPvT}.
\end{proof}

\section{Hu-Ye Conjecture on eigenvalues of hypergraphs}
We begin by introducing several fundamental classes of hypergraphs.
A \emph{hypertree} is a connected acyclic hypergraph (i.e., a connected hypergraph that contains no cycles).
For an integer $k \ge 3$ and a simple graph $G$,
the \emph{$k$-th power of $G$}, denoted by $G^{(k)}$, is a $k$-uniform hypergraph constructed as follows: each edge of $G$ (a $2$-set) is replaced by a $k$-set, with $k-2$ distinct additional vertices adjoined to each original edge to form the $k$-uniform edges of $G^{(k)}$.
Let $P_m,C_m$ and $S_m$ denote the simple path graph, the cycle graph, and the star graph with $m$ edges, respectively.
Then $P_m^{(k)}$ is called a (loose) \emph{hyperpath}, $C_m^{(k)}$ is called a \emph{hypercycle}, and $S_m^{(k)}$ is called a \emph{hyperstar}.
The \emph{complete $k$-uniform hypergraph} on $n$ vertices, denoted by $K_n^{[k]}$, is the $k$-uniform hypergraph whose edge set consists of all $k$-subsets of its vertex set.

\subsection{The spectral radius of adjacency tensor of a hypergraph}
In this section, the \emph{eigenvalues and spectral radius of a hypergraph $H$} are understood to  refer to its adjacency tensor $\A(H)$.

The subsequent results on algebraic multiplicity leverage  the Poisson formula as a foundational tool.
Combining this with the chip-firing game on graphs, Bao, Fan, et al. \cite{BaoFWZ20} derived a recursive and explicit formula for the characteristic polynomials of a class of starlike hypertrees, including hyperpaths and hyperstars.
Their work implies that the spectral radius of an $m$-edge $k$-uniform starlike hypertree has algebraic multiplicity $k^{m(k-2)}$.
Duan, van Dam, and Wang \cite{DuanvanW23} established an explicit formula of the characteristic polynomial of an $m$-edge $k$-uniform double hyperstar with $m$ edges, which shows that the
algebraic multiplicity of its spectral radius is also $k^{m(k-2)}$.
 Chen and Bu \cite{ChenB25} demonstrated a general result that unifies these findings, which is stated as follows.

\begin{thm}[\cite{ChenB25}]\label{am-tree}
The algebraic multiplicity of the spectral radius of a $k$-uniform hypertree with
$m$ edges is $k^{m(k-2)}$.
\end{thm}

Recently, Li, Su, and Fallat \cite{LSF24} expressed the characteristic polynomial of a uniform hypertree in terms of the matching polynomials of its subhypertrees, which provides a novel perspective on the algebraic multiplicity of the eigenvalues of hypertrees.

Zheng \cite{Zheng21} established an explicit formula for the characteristic polynomial of the complete $3$-uniform hypergraph $K_n^{[3]}$, which is a product of polynomials of degree at most $3$.
A key implication of his result is that the spectral radius of $K_n^{[3]}$ has algebraic multiplicity $1$ when $n >3$.

The subsequent results on algebraic multiplicity are rooted in trace-based methods.
 Using the eigenvalues of signed graphs, Duan, van Dam, and Wang \cite{DuanvanW23} established an explicit formula for the characteristic polynomial of $C_3^{(k)}$ (also called \emph{hypertriangle}), which implies that the
algebraic multiplicity of the spectral radius of $C_3^{(k)}$ is $k^{|V(C_3^{(k)})|-|E(C_3^{(k)})|-1}=k^{3k-7}$.
In a follow-up paper \cite{DuanWW24}, Duan, Wang, and Wei derived an analogous explicit formula for the characteristic polynomial of $C_4^{(k)}$, which implies that the
algebraic multiplicity of the spectral radius of $C_4^{(k)}$ is also $k^{|V(C_4^{(k)})|-|E(C_4^{(k)})|-1}=k^{4k-9}$.
Recently, Chen, van Dam, and Bu \cite{ChenvanB24} expressed the entire spectrum of $G^{(k)}$,  including multiplicities and the characteristic polynomial, in terms of parity-closed walks on the base graph $G$.
Their work unifies the aforementioned results under a general framework stated as follows.

\begin{thm}[\cite{ChenvanB24}]\label{am-power}
For a connected graph $G$,
the algebraic multiplicity of the spectral radius of  $G^{(k)}$  is
$k^{|V(G^{(k)})|-|E(G^{(k)})|-1}=k^{|V(G)|+(k-3)|E(G)|-1}.$
\end{thm}

Parts (2) and (3) of the following lemma are implicit in  \cite[Corollary 4.4]{FBH19} and \cite[Example 4.4]{FBH19}, respectively.
For completeness, we provide a full proof of the lemma by applying Theorem \ref{PV-Zk}.

\begin{lem}[\cite{FBH19}]\label{3H}
Let $H$ be a connected $k$-uniform hypergraph with spectral radius $\rho=\rho(\A(H))$. Then the following results hold.

\begin{itemize}
\item[(1)] If $H$ is a hypertree with $m$ edges, then $|\PV_\rho(\A(H))|=k^{m(k-2)}$.

\item[(2)]  If $H=G^{(k)}$ for some connected graph $G$, then $|\PV_\rho(\A(H))|=k^{|V(G)|+(k-3)|E(G)|-1}$.

\item[(3)] If $H=K_n^{[k]}$ with $n >k$, then  $|\PV_\rho(\A(H))| =1$.

\end{itemize}

\end{lem}

\begin{proof}
(1) Let $H$ be a $k$-uniform hypertree with $m$ edges.
Let $e$ be a pendent edge of $H$ and let $H-e$ be the sub-hypertree of $H$ obtained by deleting the edge $e$ together with all the vertices of $e$ that have degree $1$.
By elementary row and column transformations over $\Z_k$, the incidence matrix $B_H$ is congruent to the following direct sum: $B_H \cong 1 \oplus B_{H-e}$, where $1$ denotes a $1 \times 1$ identity matrix.
So, by induction on $m$, $B_H$ has $m$ invariant divisors all equal to $1$ over $\Z_k$.
Note that $|V(H)|=m(k-1)+1$.
By Theorem \ref{PV-Zk},
$|\PV_\rho(\A(H))|=k^{|V(H)|-m-1}=k^{m(k-2)}$.

(2) Let $H=G^{(k)}$ be the $k$-th power of a connected graph $G$.
The invariant divisors of $B_{H}$ over $\Z_k$ are $1$'s with multiplicity $|E(H)|$.
So, by Theorem \ref{PV-Zk}, $$|\PV_\rho(\A(H))|=k^{|V(H)|-|E(H)|-1}=k^{|V(G)|+|E(G)|(k-3)-1},$$ as $|V(H)|=|V(G)|+|E(G)|(k-2)$ and $|E(H)|=|E(G)|$.

(3) Let $H=K_n^{[k]}$ be the complete $k$-uniform hypergraph on $n>k$ vertices.
Consider \eqref{S0-def} for $\A(H)$:
$$\S_0(H)= \{x \in \Z_k^n: B_H x =\mathbf{0} \mbox{~over~} \Z_k, x_1=0\}.$$
It is straightforward to verify that $\S_0(H)=\{\mathbf{0}\}$, which implies the desired result by \eqref{Sadj}.
\end{proof}

\begin{proof}[\bf Proof of Theorem \ref{amPv}]
Let $\rho$ be the spectral radius of $H$.
If $H$ is a $k$-uniform hypertree or a  $k$-th power of a connected simple graph, by Theorem \ref{am-tree}, Theorem \ref{am-power}, and Lemma \ref{3H}, we have $\am(\rho)=|\PV_{\rho}(\A(H))|$.
If $H$ is a complete $3$-uniform hypergraph on at least $4$ vertices, then $\am(\rho)=1$ (\cite{Zheng21}), and hence by Lemma \ref{3H}, $\am(\rho)=|\PV_{\rho}(\A(H))|$.
For any eigenvalue $\la$ with $|\la|=\rho$, by Corollary \ref{multi-coro} and Corollary \ref{Pv-spectral}, $\am(\la)=\am(\rho)$ and $|\PV_\la(\A(H))|=|\PV_{\rho}(\A(H))|$, which implies that $\am(\la)=|\PV_{\la}(\A(H))|$.

Following the reasoning in the final part of the proof of Theorem \ref{amPvT},
the equality in \eqref{conj1} is directly valid for all eigenvalues $\la$ of $\A(H)$ with the modulus $\rho$.
The conclusion thus follows.
\end{proof}

\subsection{The zero eigenvalue of the Laplacian or signless Laplacian tensor of a hypergraph}
Recently, Zheng \cite{Zheng24-2} proved that the algebraic multiplicity of the zero eigenvalue of the Laplacian tensor $\L(H)$ is $k^{m(k-2)}$ for any $k$-uniform hyperpath or hyperstar $H$ with $m$ edges.
Lin and Bu \cite{LinB23} further generalized this result to arbitrary $k$-uniform hypertrees.

\begin{thm}[\cite{LinB23}]\label{zeroL}
Let $T$ be a $k$-uniform hypertree with $m$ edges.
Then the algebraic multiplicity of the zero eigenvalue of the Laplacian tensor $\L(T)$ is $k^{m(k-2)}$.
\end{thm}

\begin{proof}[\bf Proof of Theorem \ref{LamPv}]
Let $T$ be a $k$-uniform hypertree.
By Theorem \ref{am-tree} and Theorem \ref{zeroL}, we have
\begin{equation}\label{zerorho} \am(0, \L(T))=\am(\rho(\A(T)), \A(T)),
\end{equation}
where $\rho(\A(T))$ is abbreviated as $\rho$ for conciseness in prior sections.
For any connected $k$-uniform hypergraph $H$, combining \eqref{Sadj} and \eqref{Slap}, we have
\begin{equation}\label{Vzerorho} |\PV_0(\L(H))|=|\PV_\rho(\A(H))|.
\end{equation}
So, the result follows by Theorem \ref{amPv} immediately.
\end{proof}

We conjecture that the equality \eqref{zerorho} extends to all connected $k$-uniform hypergraphs, not just hypertrees.

\begin{conj}\label{conjzerorho}
Let $H$ be a connected $k$-uniform hypergraph. Then
\begin{equation}\label{conjzero1} \am(0, \L(H))= \am(\rho(\A(H)), \A(H)).
\end{equation}
\end{conj}

Suppose that Conjecture \ref{conjzerorho} and Conjecture \ref{conjH} hold.
Combining these with Eq. \eqref{Vzerorho}, we will have
\begin{equation}\label{conjzero2} \am(0, \L(H))= |\PV_0(\L(H))|.
\end{equation}

We conclude with a discussion of the signless Laplacian tensor of a connected $k$-uniform hypergraph $H$.
Fan et al. \cite{FanW19} proved that $\Q(H)$ has a zero eigenvalue if and only if there exists a nonsingular diagonal matrix $D$ such that $\Q(H)=D^{-(k-1)}\L(H)D$.
In this case, $\Q(H)$ and $\L(H)$ are similar and hence have identical spectra \cite{Shao13}, which implies that
$$ \am(0,\Q(H))=\am(0,\L(H)).$$
By Eq. \eqref{Sslap}, $|\PV_0(\Q(H))|=|\PV_0(\L(H))|$ whenever $\Q(H)$ has a zero eigenvalue.
Thus, if Eq. \eqref{conjzero2} holds, we will have
$$ \am(0,\Q(H))=|\PV_0(\Q(H))|.$$

\section{Simple zeros of polynomials}
\subsection{Local intersection multiplicity}\label{MultiZero}
Throughout this section, let $\K=\C$.
For a point $p=(p_1,\ldots,p_n)\in \K^n$, denote $M_p=\langle x_1-p_1,\ldots,x_{n}-p_{n}\rangle$, the maximal ideal in $\K[x_1,\ldots,x_{n}]$ corresponding to $p$.
The \emph{localization of $\K[x_1,\ldots,x_n]$ at $M_p$} is defined as
$$\K[x_1,\ldots,x_{n}]_{M_p}=\{f/g: f,g \in \K[x_1,\ldots,x_n], g(p) \ne 0\},$$
where $\K(x_1,\ldots,x_n)$ denotes the field of rational functions over $\K$.
It is well-known that $\K[x_1,\ldots,x_{n}]_{M_p}$ is the so-called \emph{local ring}.
Let $J$ be a zero-dimensional ideal in $\K[x_1,\ldots,x_n]$ (i.e., its zero set $\V(J)\subseteq\K^n$ is finite), and write
\[
 J_p:=J\K[x_1,\ldots,x_n]_{M_p}
\]
for the ideal of the local ring generated by $J$.
For  $p \in \V(J)$,
the \emph{local intersection multiplicity} of $p$ with respect to $J$ is defined as
$$ m(p;J)=\dim_{\K} \K[x_1,\ldots,x_{n}]_{M_p}/J_p,$$
where the right-hand side denotes the dimension of the quotient ring as a $\K$-vector space.
When the ideal $J$ is clear from the context, we abbreviate this notation to $m(p)$ for conciseness.

An equivalent characterization of the local intersection multiplicity can be given using differential functionals (dual space formulation).
For an index array $j=(j_1,\ldots,j_n)\in \N^n$ and $x=(x_1,\ldots,x_n)$, let
$x^j=x_1^{j_1}\cdots x_n^{j_n}$ and $(x-y)^j=(x_1-y_1)^{j_1}\cdots (x_n-y_n)^{j_n}$.
The \emph{order} of $j$ is $|j|=j_1+\cdots+j_n$.
Define the \emph{normalized partial differential operator} of type $j$ as
$$\part_j=\part_{j_1\ldots j_n}=\frac{1}{j_1!\ldots j_n!}\frac{\part^{j_1+\cdots+j_n}}{\part x_1^{j_1} \cdots \part x_n^{j_n}}.$$
For a point $p \in \K^n$, the \emph{differential functional at $p$ of type $j$} is the linear function:
$$ \part_j[p]: \C[x_1,\ldots,x_n] \to \C, f \mapsto (\part_j f)(p).$$
The \emph{dual space of $J$ at $p$} is defined as
$$\mathcal{D}_p(J)=\{c=\sum_{j \in \N^n} c_j \part_j[p]: c(f)=0 \mbox{~for all ~} f \in J\}.$$
For $\alpha \in \N$, let $\mathcal{D}^\alpha_p(J) \subseteq \mathcal{D}_p(J)$ denote the subset of functionals with differential orders bounded by $\alpha$.
The \emph{Hilbert function} of $\mathcal{D}_p(J)$ is defined by
$$ H(0)=\dim \mathcal{D}^0_p(J)=1, H(\alpha)=\dim \mathcal{D}^\alpha_p(J) - \dim \mathcal{D}^{\alpha-1}_p(J) \mbox{~for~} \alpha \ge 1.$$
Let $\sigma \in \N$ be the smallest integer such that $\dim_\K \mathcal{D}^\alpha_p(J) = \dim_\K \mathcal{D}^{\alpha+1}_p(J)$ (or equivalently, $H(\alpha+1)=0)$ for all $\alpha \ge \sigma$).
Then $\mathcal{D}_p(J)=\mathcal{D}^\sigma_p(J)$.

Dayton and Zeng \cite{DayZ05} proved that for a zero-dimensional ideal $J \subseteq \K[x_1, \ldots, x_n]$ and a point $p \in \V(J)$, the local intersection multiplicity $m(p;J) =m$ if and only if $\dim_\K \mathcal{D}_p(J)=m$.
Furthermore,  the authors introduced multiplicity matrices and employed these structural tools to characterize Hilbert functions.

\begin{thm}[\cite{DayZ05}]\label{mul-mat}
Let $p$ be an isolated zero of $f_1,\ldots,f_t \in \K[x_1,\ldots,x_n]$, where $t \ge n$.
For $\alpha \in \N$, a differential functional $\sum_{|j| \le \alpha} c_j \part_j[p]$ belongs to
$\mathcal{D}_p^\alpha(\langle f_1,\ldots, f_t\rangle)$ if and only if the coefficient vector $c=[c_j: |j| \le \alpha] \in \K^{n_\alpha}$ satisfies
$$ \sum_{|j| \le \alpha} c_j \part_j((x-p)^k f_i)(p)=0, k \in \N^n, |k|<\alpha, i \in [t],$$
which corresponds to a homogeneous linear system $S_\alpha c=0$, where $S_\alpha$ is an $m_\alpha \times n_\alpha$ matrix with rows indexed by $(k,i)$ and columns indexed by $j$, and is called the $\alpha$-th order multiplicity matrix.
Here $m_\alpha=\binom{\alpha-1+n}{\alpha -1}\cdot t$ and $n_\alpha=\binom{\alpha+n}{\alpha}$.
Consequently, the Hilbert function $$H(\alpha)=\mbox{nullity}(S_\alpha)-\mbox{nullity}(S_{\alpha-1})$$ for $\alpha \ge 1$, with $S_0=[f_1(p),\ldots,f_t(p)]^\top=O_{t \times 1}$ (the $t \times 1$ zero matrix).
\end{thm}

A direct consequence of Theorem \ref{mul-mat} is the following characterization of \emph{simple zeros} (i.e. isolated zeros with local intersection multiplicity $1$).

\begin{coro}\label{simple}
Let $p \in \K^n$ be an isolated zero of $f_1,\ldots,f_t \in \K[x_1,\ldots,x_n]$, where $t \ge n$.
Then $p$ is a simple zero if and only if the Jacobian matrix
$$J(f_1,\ldots,f_t)=\left(\frac{\part f_i}{\part x_j}\right)_{t \times n}$$
has rank $n$ at $p$.
\end{coro}

\begin{proof}
The multiplicity matrix $S_1$ in Theorem \ref{mul-mat} can be written as
$$ S_1=\left[\begin{array}{c|c}
f_1(p) &    \\
\vdots &   J(p)\\
f_t(p)  &
\end{array}
\right]=
\left[\begin{array}{c|c}
0 &    \\
\vdots &   J(p)\\
0  &
\end{array}
\right],
$$
where $J(p)=J(f_1,\ldots,f_t)(p)$ denotes the Jacobian matrix evaluated at $p$.
If $p$ is simple, then $m(p)=1$.  Theorem \ref{mul-mat} and the definition of the Hilbert function give $H(0)=1$ and $H(\alpha)=0$ for every $\alpha\ge1$.  In particular,
$H(1)=\mbox{nullity}(S_1)-\mbox{nullity}(S_0)=0$, so $\mbox{nullity}(S_1)=1$, which is equivalent to $J(p)$ having rank $n$.

Conversely, if $J(p)$ has rank $n$, then
\[
H(1)=\operatorname{nullity}(S_1)
-\operatorname{nullity}(S_0)=0.
\]
By \cite[Lemma~1(iv)]{DayZ05}, $H(1)=0$ implies that
$H(\alpha)=0$ for all $\alpha\geq1$.
Therefore,
\[
m(p)=\dim_{\K}\mathcal D_p(J)
=\sum_{\alpha\geq0}H(\alpha)=H(0)=1.
\]
Hence $p$ is a simple zero.
\end{proof}

We note that Corollary \ref{simple} is equivalent to the Jacobian criterion for simple zeros (cf. \cite[Theorem 16.19]{Eisen}).
Our proof leverages Theorem \ref{mul-mat} to avoid relying on additional algebraic geometry prerequisites.

\subsection{Simple zeros of projective eigenvariety}
For an isolated point $p$ of a projective variety defined by a
homogeneous ideal, its local intersection multiplicity is defined
as the local intersection multiplicity of the corresponding
dehomogenized ideal on any affine chart containing $p$.
This number is independent of the chosen affine chart.
In this section, we establish a sufficient condition for an isolated point in a projective eigenvariety of a weakly irreducible symmetric tensor to have local intersection multiplicity $1$.

\begin{lem}\label{simple2}
Let $\A=(a_{i_1\ldots i_k})$ be a weakly irreducible symmetric tensor of order $k\ge2$ and dimension $n\ge2$ such that $a_{i_1\ldots i_k}\ne0$ only when the indices $i_1,\ldots,i_k$ are either all equal or pairwise distinct.  For every $k$-set $e=\{i_1,\ldots,i_k\}$ corresponding to a nonzero off-diagonal entry, define
\[
 c_e=(k-1)!\,a_{i_1\ldots i_k};
\]
this is well defined because $\A$ is symmetric.  Let $\la$ be an eigenvalue of $\A$, and let $p\in\PV_\la(\A)$ be an isolated point.  Suppose that
\begin{itemize}
\item[(1)] $p$ has no zero entries;
\item[(2)] there is an $\alpha\in\mathbb{R}$ such that
\[
 c_ep^e=|c_ep^e|e^{\i\alpha}
\]
for every such $k$-set $e$, where $p^e=\prod_{i\in e}p_i$.
\end{itemize}
Then the local intersection multiplicity of $p$ with respect to the homogeneous eigen-equation ideal defining $\PV_\la(\A)$ is $1$.
\end{lem}

\begin{proof}
Let $H$ be the $k$-uniform hypergraph on $[n]$ whose edge set is
\[
 E=\{e\subseteq[n]:|e|=k,\ c_e\ne0\}.
\]
The symmetry and weak irreducibility of $\A$ imply that $H$ is connected.  For $v\in[n]$, let $E_v$ be the set of edges containing $v$, and for distinct $v,w\in[n]$, let $E_{v,w}$ be the set of edges containing both vertices.  The eigenvector equations are
\begin{equation}\label{weighted-Fi}
 F_v(x)=(\la-a_{v,\ldots,v})x_v^{k-1}
 -\sum_{e\in E_v}c_ex^{e\setminus\{v\}}=0,
 \quad v\in[n].
\end{equation}
Indeed, the factor $(k-1)!$ in $c_e$ accounts for the permutations of the remaining $k-1$ indices in $(\A x^{k-1})_v$.

By Condition (1), normalize $p_n=1$; Condition (2) is preserved, with at most a common change in $\alpha$, because every edge has size $k$.  Put $\tilde p=(p_1,\ldots,p_{n-1})$, and set
\[
 f_v(x_1,\ldots,x_{n-1})=F_v(x_1,\ldots,x_{n-1},1),
 \quad v\in[n].
\]
Then $\tilde p$ is a common zero of $f_1,\ldots,f_n$.  The entries of the Jacobian
$J(\tilde p)=J(f_1,\ldots,f_{n-1})(\tilde p)$ are
\begin{equation}\label{weighted-Jacobian}
 J(\tilde p)_{vv}=(k-1)(\la-a_{v,\ldots,v})p_v^{k-2},\quad
 J(\tilde p)_{vw}=-\sum_{e\in E_{v,w}}c_ep^{e\setminus\{v,w\}}\quad(v\ne w),
\end{equation}
where $p_n=1$ is understood.

Let $G_n$ be the graph on $[n-1]$ in which distinct vertices $v,w$ are adjacent if and only if $E_{v,w}\ne\emptyset$, and let $C_1,\ldots,C_s$ be its connected components.  After ordering the vertices component by component, $J(\tilde p)$ is block diagonal, with a block $J_i(\tilde p)$ indexed by $V(C_i)$.  For distinct $v,w\in V(C_i)$, Condition (2) gives
\begin{equation}\label{weighted-offdiag}
 J_i(\tilde p)_{vw}
 =-\frac{1}{p_vp_w}\sum_{e\in E_{v,w}}c_ep^e
 =-\frac{e^{\i\alpha}}{p_vp_w}
   \sum_{e\in E_{v,w}}|c_ep^e|.
\end{equation}
Thus $J_i(\tilde p)_{vw}\ne0$ exactly when $E_{v,w}\ne\emptyset$.  The graph of its nonzero off-diagonal entries is exactly $C_i$; hence $J_i(\tilde p)$ is irreducible.

Fix a component $C_i$ and put $D_i=\diag(p_v:v\in V(C_i))$.  This diagonal matrix is nonsingular.  For $v\in V(C_i)$, by \eqref{weighted-offdiag}, double counting, and \eqref{weighted-Fi}, we have
\begin{equation}\label{weighted-dominance}
\begin{split}
\sum_{w\in V(C_i)\setminus\{v\}}|(J_i(\tilde p)D_i)_{vw}|
&=\frac{1}{|p_v|}\sum_{e\in E_v}|c_ep^e|\,|e\setminus\{v,n\}| \\
&\le\frac{k-1}{|p_v|}\sum_{e\in E_v}|c_ep^e| \\
&=\frac{k-1}{|p_v|}\left|\sum_{e\in E_v}c_ep^e\right| \\
&=(k-1)\left|\sum_{e\in E_v}c_ep^{e\setminus\{v\}}\right| \\
&=(k-1)|\la-a_{v,\ldots,v}|\,|p_v|^{k-1}
=|(J_i(\tilde p)D_i)_{vv}|.
\end{split}
\end{equation}
Thus $J_i(\tilde p)D_i$ is diagonally dominant.

Every component $C_i$ of $G_n$ is joined to vertex $n$ in the two-section of $H$ because $H$ is connected.  Hence there is an edge $e_0\in E$ such that
\[
 n\in e_0,\quad e_0\setminus\{n\}\subseteq V(C_i).
\]
Choose $v_0\in e_0\setminus\{n\}$.  
Since $|c_{e_0}p^{e_0}|>0$ and
$|e_0\setminus\{v_0,n\}|=k-2<k-1$, the inequality in
\eqref{weighted-dominance} is strict for $v=v_0$.
Consequently, the $v_0$-th row is strictly diagonally dominant.
Since $J_i(\tilde p)D_i$ is also irreducible, it is nonsingular by \cite[Corollary 6.2.9]{HornJ}.  Therefore every $J_i(\tilde p)$ is nonsingular, and $J(\tilde p)$ has rank $n-1$.
Hence the full Jacobian
$J(f_1,\ldots,f_n)(\tilde p)$ also has rank $n-1$,
since it contains $J(\tilde p)$ as a nonsingular
$(n-1)\times(n-1)$ submatrix.

Since $p$ is isolated in $\PV_\la(\A)$, its dehomogenization
$\tilde p$ is an isolated common zero of $f_1,\ldots,f_n$.
By Corollary \ref{simple}, $\tilde p$ is a simple zero of this
polynomial system. Since dehomogenization on the chart $x_n\ne0$
preserves local intersection multiplicities, the local intersection
multiplicity of $p$ with respect to the homogeneous eigen-equation
ideal is $1$.
\end{proof}

\begin{proof}[\bf Proof of Theorem \ref{multirho}]
To establish the result, it suffices to verify that the adjacency tensor $\A(H)$ and its eigenvectors satisfy the hypotheses of Lemma \ref{simple2}.
Since $H$ is a connected $k$-uniform hypergraph, $\A(H)$ is symmetric and weakly irreducible, and every edge weight in Lemma \ref{simple2} is $c_e=1$.
Let $\la$ be an eigenvalue of $\A(H)$ with $|\la|=\rho(\A(H))$.
By Theorem \ref{PFthm}, each point $p \in \PV_\la(\A(H))$ has no zero entries, satisfying Condition (1) of Lemma \ref{simple2}.
By the uniform angle property (cf. \eqref{angle}), the products $c_ep^e=p^e$ have a common argument, so Condition (2) of Lemma \ref{simple2} holds.
By results in \cite{FanHB22} or Corollary \ref{Pv-spectral}, $\PV_{\la}(\A(H))$ is finite.
So $p$ is an isolated point in $\PV_{\la}(\A(H))$.
The result now follows.
\end{proof}

Similarly, we establish analogous multiplicity results for the signless Laplacian tensor $\Q(H)$ and the Laplacian tensor $\L(H)$ of a connected hypergraph $H$.
By Theorem \ref{symequ}, the cyclic index of $\Q(H)$ is equal to $1$  as $\Tr_1(\Q(H)) >0$, which implies that its spectral radius $\rho(\Q(H))$ is the unique eigenvalue of $\Q(H)$ with modulus $\rho(\Q(H))$.

\begin{coro}\label{multirho-Q}
Let $H$ be a connected $k$-uniform hypergraph, and let $\rho$ be the spectral radius of $\Q(H)$.
Then each point $p \in \PV_{\rho}(\Q(H))$ has local intersection multiplicity $1$.
\end{coro}

\begin{proof}
We verify the hypotheses of Lemma \ref{simple2} for $\Q(H)$ and its eigenvectors corresponding to $\rho$.
 Since $H$ is connected, $\Q(H)$ is symmetric and weakly irreducible, with edge weights $c_e=1$.
 By Theorem \ref{PFthm}, each point $p \in \PV_{\rho}(\Q(H))$ has no zero entries; and by \eqref{angle}, the products $c_ep^e=p^e$ have a common argument.
 By the  results in \cite{FanHB22} or Corollary \ref{Pv-spectral}, $p$ is an isolated point.
 All hypotheses of Lemma \ref{simple2} are satisfied.
 Hence the local intersection multiplicity of $p$ is $1$.
\end{proof}

Next, we consider the Laplacian tensor $\L(H)$ of $H$.

\begin{coro}
Let $H$ be a connected $k$-uniform hypergraph.
Then each point $p \in \PV_{0}(\L(H))$ has local intersection multiplicity $1$.
\end{coro}

\begin{proof}
For $\L(H)$, the edge weights in Lemma \ref{simple2} are $c_e=-1$.  Equation \eqref{Lzero} shows that $p$ has no zero entries, while \eqref{unimod} gives $p^e=1$ for every edge $e$.  Thus $c_ep^e=-1$ has the same argument for all edges.  The finiteness of $\PV_0(\L(H))$ in \eqref{Slap} shows that $p$ is isolated.  Hence Lemma \ref{simple2} applies, and the local intersection multiplicity of $p$ is $1$.
\end{proof}

\begin{coro}
Let $H$ be a connected $k$-uniform hypergraph.
If $\Q(H)$ has a zero eigenvalue, then each point $p \in \PV_{0}(\Q(H))$ has local intersection multiplicity $1$.
\end{coro}

\begin{proof}
For $\Q(H)$, the edge weights in Lemma \ref{simple2} are $c_e=1$.  Equation \eqref{Lzero} shows that $p$ has no zero entries, and \eqref{Qunimod} gives $p^e=-1$ for every edge $e$.  Hence the products $c_ep^e=-1$ have a common argument.  Finally, \eqref{Sslap} implies that $p$ is isolated in $\PV_0(\Q(H))$.  The result follows from Lemma \ref{simple2}.
\end{proof}

\subsection{A Smith-form reformulation of Conjecture \ref{conjH}}
Let $H$ be a connected $k$-uniform hypergraph on $n$ vertices, set $\A=\A(H)$ and $\rho=\rho(\A)$, and retain the notation $t,R,Q,M,M_0$, and $Z_\rho$ from the proof of Theorem \ref{amPvT}.  Put $N=(k-1)^{n-1}$.  The Smith normal form of $M$ over the discrete valuation ring $R$ can be written as
\begin{equation}\label{smith-M}
 UMV=\diag(t^{r_1},\ldots,t^{r_s},1,\ldots,1),
 \quad r_i\ge1,
\end{equation}
where $U,V\in\operatorname{GL}_N(R)$.  It follows that
\begin{equation}\label{smith-orders}
 \am(\rho)=\sum_{i=1}^s r_i,
 \quad \dim_{\C}\ker M_0=s.
\end{equation}

We next identify this nullity geometrically.  The Artinian algebra $Q/tQ$ is the product of its local factors.  On a factor corresponding to a point outside $Z_\rho$, multiplication by $f_n|_{\la=\rho}$ is invertible.  If $q\in Z_\rho$, let $\mathfrak m_q$ be the corresponding maximal ideal.  The dimension of the kernel on the local factor equals the dimension of its cokernel, and the latter local algebra is
\[
 \frac{(Q/tQ)_q}{f_n|_{\la=\rho}(Q/tQ)_q}
 \cong
 \frac{\C[x_1,\ldots,x_{n-1}]_{M_q}}
 {\langle f_1|_{\la=\rho},\ldots,f_n|_{\la=\rho}\rangle
  \C[x_1,\ldots,x_{n-1}]_{M_q}}.
\]
Under the bijection $Z_\rho\cong\PV_\rho(\A)$ in the proof of Theorem \ref{amPvT}, the dimension of this algebra is the local intersection multiplicity of the corresponding projective eigenvector.  By Theorem \ref{multirho}, that multiplicity is $1$.  Therefore,
\begin{equation}\label{smith-nullity}
 \dim_{\C}\ker M_0=|Z_\rho|=|\PV_\rho(\A)|.
\end{equation}
Combining \eqref{smith-orders} and \eqref{smith-nullity}, we obtain the precise criterion
\[
 \am(\rho)=|\PV_\rho(\A)|
 \quad\Longleftrightarrow\quad
 r_1=\cdots=r_s=1.
\]
Thus Conjecture \ref{conjH} is equivalent to the assertion that multiplication by the last dehomogenized eigen-equation has no Smith elementary divisor divisible by $t^2$.  This formulation permits collisions in the special fiber and does not require specializing individual algebraic branches.

\end{document}